\documentclass{amsart}

\usepackage{amsmath,amssymb,amsthm}

\title{Who owns the theorem?}

\author{Melvyn B. Nathanson}\address{Lehman College (CUNY), Bronx, NY 10468}\email{melvyn.nathanson@lehman.cuny.edu}

\subjclass[2010]{00A30,01A80.}

\keywords{Philosophy of mathematics, sociology of mathematics}

\date{\today}
\begin{document}

\begin{abstract}
Epistemological and sociological questions about ownership of mathematical and scientific discoveries.
\end{abstract}

\maketitle

Simple questions: 
If you prove the theorem, do you own it?   
Can you forbid others to use or even cite it?  
Can you choose not publish the theorem?  
Can you be forbidden to publish it?  

What is theorem?   
A mathematical statement may be true.  It is true whether or not there is a proof.  
Without a proof, we do not know if it is true.  
A theorem is a true mathematical statement that has a proof.  

Suppose there is a true mathematical statement, and you prove it.  Now it is a theorem.  
It is ``your theorem.''   In what sense might you own it?  
Can or should a theorem be considered the private property of its discoverer, 
who may or may not choose to publish?  
If you own the theorem, can you license it or rent it?  
Can you insist that anyone who wants to use or apply the theorem must pay you to do so?

If you publish the theorem in a refereed journal, or post it on arXiv,  
or explain it in a seminar, or submit it to a journal, then everyone knows you proved it.  
When does ``your theorem'' become part of the public library of proven mathematical 
truths that other researchers can freely use to prove new theorems?  

If you need a result to prove a theorem, 
and know that the result is true but the discoverer has not announced or 
released it publicly, is it ethical (of course, properly citing the discoverer) 
to use that ``unpublished'' result in the  proof?
Is it ethical for you \textit{not} to prove the theorem because it requires a result 
that is true but is being withheld by its ``owner''?

Suppose you find out that someone has proved a theorem, but has not revealed it to the world.
Maybe you have even seen the proof, and checked it, so you are sure that it is correct.
Even though it has not been published, you know that it is a mathematical truth.  
Can you use it in a paper, even though the discoverer might not 
want the result to be known?  Does the prover of the theorem own it enough to prevent 
other mathematicians from using it?  

The notion of ``owning a mathematical truth'' is, in part, connected with careerism in academic life.  
What might be called ``vulgar careerism''  is endemic and  not necessarily vulgar.  
Many mathematicians hide what they are working on so others will not ``scoop'' them, will not 
use ``their'' ideas to prove a theorem before they do.  
Perhaps it is not sufficient to give proper attribution.  
Maybe the author is an untenured assistant professor who wants to deduce more results from the theorem,  
publish more papers, and get promoted.  
Maybe the author thinks it will lead to a proof of the Riemann hypothesis 
and earn the million dollar prize from the Clay Mathematics Institute.  
Some mathematicians admit that they discuss their ideas about how to solve the Riemann hypothesis 
only after they are convinced that the ideas will not work.  

It used to be that, every year, permanent professors in mathematics at the Institute for Advanced Study 
in Princeton would appoint a visiting member to be their ``assistant.''  
Long ago at the Institute, there was a permanent member 
who required his assistants to promise not to reveal to anyone what he was working on.\footnote{I was once 
Andr\' e Weil's assistant at the Institute.  He did not impose a secrecy oath.} 
When I learned this, I was shocked.  It was antithetical to everything I believed about science.  I was also naive.  
 I had not understood that for many people mathematics is a competition.\footnote{I still do not understand why, 
 for some mathematicians, getting medals in high school and college competitions is a core part of their self-esteem.} 

In 1977,  the National Security Agency decided that publication of cryptographic research 
would endanger national security,  and wanted to require that 
professors who wrote papers in cryptography would have to send them 
for pre-publication review by the NSA and not submit them to journals without NSA approval. 
At first, NSA hoped for voluntary compliance, but also considered making this a legal requirement.  
This did not happen, and, after considerable contentiousness and debate 
in the mathematical community,  prepublication review  by the NSA faded away and 
seems not to be an issue today~\cite{land83,diff-land07}.

Secrecy in mathematics is less important than in other sciences.
 Mathematical results rarely have commercial value.  
 Like many mathematicians, I don't care if my theorems are  ``useful.''  
 I only hope that I have not made mistakes, that the proofs are correct, 
 that the ``theorems'' are theorems and are interesting.  
 I upload preprints to arXiv as soon as they are written, 
 before I send them to a journal.  I am happy if someone uses my results.  
 But this does not answer the central question.  
 
Mathematician A proves a theorem, and mathematician B learns about it.  
Maybe B reviewed A's NSF proposal. 
Maybe A submitted the manuscript to  a journal and B refereed it.  
Can B use the theorem (as always, with proper attribution) in a paper before A has published it? 

For me, the answer is clear.  
Here is an analogy.  
Legally, you cannot sequence a plant or animal DNA strand and patent the sequence 
because you did not create the sequence.  God created it, and you only discovered it.  Similarly, 
mathematical truths exist, and mathematicians only discover them.  If you discover a theorem, you have the power, the privilege, and, perhaps, the right not to reveal it to anyone, 
but if, somehow, someone learns of your result, 
knows that a certain mathematical statement is true, then that person has the right 
to tell the world and to apply it to obtain new results, with or without your consent.  

Can you own a scientific  truth? Can you hide a scientific truth? 
These are ethical questions, and, in the covid era, not only in mathematics.  

\def\cprime{$'$} \def\cprime{$'$}
\providecommand{\bysame}{\leavevmode\hbox to3em{\hrulefill}\thinspace}
\providecommand{\MR}{\relax\ifhmode\unskip\space\fi MR }
\providecommand{\MRhref}[2]{%
  \href{http://www.ams.org/mathscinet-getitem?mr=#1}{#2}
}
\providecommand{\href}[2]{#2}

\end{document}